# Improved estimation of population variance using information on auxiliary attribute in simple random sampling

Rajesh Singh and Sachin Malik

Department of Statistics, Banaras Hindu University

Varanasi-221005, India

(rsinghstat@gmail.com, sachinkurava999@gmail.com)

**Abstract**

Singh and Kumar (2011) suggested estimators for calculating population variance using auxiliary attributes. This paper proposes a family of estimators based on an adaptation of the estimators presented by Kadilar and Cingi (2004) and Singh et al. (2007), and introduces a new family of estimators using auxiliary attributes. The expressions of the mean square errors (MSEs) of the adapted and proposed families are derived. It is shown that adapted estimators and suggested estimators are more efficient than Singh and Kumar (2011) estimators. The theoretical findings are supported by a numerical example.

## 1. Introduction

It is well known that the auxiliary information in the theory of sampling is used to increase the efficiency of estimator of population parameters. Out of many ratio, regression and product methods of estimation are good examples in this context. There exist situations when information is available in the form of attribute which is highly correlated with y. Taking into consideration the point bi-serial correlation coefficient between auxiliary attribute and study variable, several authors including Naik and Gupta (1996), Jhajj et al. (2006), Shabbir and Gupta (2007), Singh et al. (2007, 2008), Singh et al. (2010), Abd-Elfattah et al. (2010), Malik and Singh (2013a,b,c), Singh (2013), Singh and Malik (2013) and Sharma et al. (2013) defined ratio estimators of population mean when the prior information of population proportion of units, possessing the same attribute is available.

In many situations, the problem of estimating the population variance $\sigma^2$ of study variable y assumes importance. When the prior information on parameters of auxiliary variable(s) is available, Das and Tripathi (1978), Isaki (1983), Prasad and Singh (1990), Kadilar and Cingi (2006, 2007) and Singh et al. (2007) have suggested various estimators of $S_y^2$.

Consider a sample of size n drawn by SRSWOR from a population of size N. Let $y_i$ and $\varphi_i$ denote the observations on variable y and $\varphi$ respectively for the $i^{th}$ unit (i=1,2,3....N). It is assumed that attribute $\varphi$ takes only the two values 0 and 1 according as

$\varphi = 1$, if $i^{th}$ unit of the population possesses attribute $\varphi$

$= 0$, if otherwise.

The variance of the usual unbiased estimator $\hat{S}_y^2$ is given by

$$\text{var}(\hat{S}_y^2) = \frac{S_y^4}{n}(\lambda_{40} - 1) \tag{1.1}$$

where, $\lambda_{rq} = \frac{\mu_{rq}}{\left(\mu_{20}^{r/2}\mu_{02}^{q/2}\right)}$, $\mu_{rq} = \frac{\sum_{i=1}^{N}(y_i - \overline{Y})^r(\varphi_i - P)^q}{(N-1)}$.

In this paper we have proposed a family of estimators for the population variance $S_y^2$ when auxiliary variable is in the form of attribute. For main results we confine ourselves to sampling scheme SRSWOR ignoring the finite population correction.

## 2. Estimators in literature

In order to have an estimate of the study variable y, assuming the knowledge of the population proportion P, Singh and Kumar (2011) proposed the following estimators

$$t_1 = s_y^2 \frac{S_\phi^2}{s_\phi^2} \tag{2.1}$$

$$t_2 = s_y^2 + b_\varphi\left(S_\varphi^2 - s_\varphi^2\right) \tag{2.2}$$

$$t_3 = s_y^2 \exp\left[\frac{S_\varphi^2 - s_\varphi^2}{S_\varphi^2 - s_\varphi^2}\right] \quad (2.3)$$

The MSE expression of the estimator $t_1$ and variance of $t_2$ are given, respectively, by

$$MSE(t_1) = \frac{S_y^4\left[(\lambda_{40} - 1) + (\lambda_{04} - 1) - 2(\lambda_{22} - 1)\right]}{n} \quad (2.4)$$

$$V(t_2) = \frac{1}{n}\left[S_y^4(\lambda_{40} - 1) + b_\phi^2 S_\varphi^2(\lambda_{04} - 1) - 2b_\phi S_y^2 S_\varphi^2(\lambda_{22} - 1)\right] \quad (2.5)$$

On differentiating (2.5) with respect to b and equating to zero we obtain

$$b_\varphi = \frac{S_y^2(\lambda_{22} - 1)}{S_x^2(\lambda_{04} - 1)} \quad (2.6)$$

Substituting the optimum value of $b_\varphi$ in (2.5), we get the minimum variance of the estimator $t_2$, as

$$\min. V(t_2) = \frac{S_y^4}{n}\left[(\lambda_{40} - 1) - \frac{(\lambda_{22} - 1)^2}{(\lambda_{04} - 1) - 1}\right] \quad (2.7)$$

The MSE expression of the estimator $t_3$ is given by

$$MSE(t_3) = \frac{S_y^4}{n}\left[(\lambda_{40} - 1) + \frac{(\lambda_{04} - 1)}{4} - (\lambda_{22} - 1)\right] \quad (2.8)$$

## 3. The adapted estimators

Following Kadilar and Cingi (2004), we propose the following variance estimators using known values of some population parameter(s),

$$t_{KC1} = s_y^2\left(\frac{S_\varphi^2 + C_p}{s_\varphi^2 + C_p}\right) \quad (3.1)$$

$$t_{KC2} = s_y^2\left(\frac{S_\varphi^2 + \beta_{2\varphi}}{s_\varphi^2 + \beta_{2\varphi}}\right) \quad (3.2)$$

$$t_{KC3} = s_y^2 \left( \frac{S_\varphi^2 \beta_{2\varphi} + C_p}{s_\varphi^2 \beta_{2\varphi} + C_p} \right) \tag{3.3}$$

$$t_{KC4} = s_y^2 \left( \frac{S_\varphi^2 C_p + \beta_{2\varphi}}{s_\varphi^2 C_p + \beta_{2\varphi}} \right) \tag{3.4}$$

where, $s_y^2$ and $s_\phi^2$ are unbiased estimator of population variances $S_y^2$ and $S_\phi^2$ respectively.

To obtain the bias and MSE, we write-

$$s_y^2 = S_y^2(1+e_0), \qquad s_\phi^2 = S_\phi^2(1+e_1)$$

Such that $E(e_0) = E(e_1) = 0$

and $E(e_0^2) = \dfrac{(\lambda_{40}-1)}{n}$, $E(e_1^2) = \dfrac{(\lambda_{04}-1)}{n}$, $E(e_0 e_1) = \dfrac{(\lambda_{22}-1)}{n}$,

and $k_{pb} = \rho_{pb} \dfrac{C_y}{C_p}$.

The MSE expressions of $\hat{S}_{KC_i}$ $(i=1,2,3,4)$ to the first order of approximation are respectively given by

$$\text{MSE}(t_{KCi}) = \frac{S_y^4}{n}\left[(\lambda_{40}-1) + \omega_i^2(\lambda_{04}-1) - 2\omega_i(\lambda_{22}-1)\right], (i=1,2,3,4) \tag{3.5}$$

where, $\omega_1 = \dfrac{S_\varphi^2}{S_\varphi^2 + C_p}$, $\omega_2 = \dfrac{S_\varphi^2}{S_\varphi^2 + \beta_{2\varphi}}$, $\omega_3 = \dfrac{S_\varphi^2 \beta_{2\varphi}}{S_\varphi^2 \beta_{2\varphi} + C_p}$, $\omega_4 = \dfrac{S_\varphi^2 C_p}{S_\varphi^2 C_p + \beta_{2\varphi}}$.

Following Singh et al. (2007), we propose estimator $t_S$ as

$$t_S = \frac{s_y^2 + b_\varphi(S_\varphi^2 - s_\varphi^2)}{(n_1 s_\varphi^2 + n_2)}(n_1 S_\varphi^2 + n_2) \tag{3.6}$$

where $n_1, n_2$ are either real numbers or the functions of the known parameters of attribute such as $C_p, \beta_{2\varphi}, \rho_{pb}$ and $k_{pb}$.

Expressing equation (3.6) in terms of e's, we have

$$t_s = \left(S_y^2(1+e_0) - b_\varphi S_\varphi^2 e_1\right)\left[1 + A_1 e_1^2\right]^{-1}$$

where, $A_1 = \dfrac{n_1 S_\varphi^2}{n_1 S_\varphi^2 + n_2}$.

Up to first order of approximation, the MSE of $t_s$ is given by

$$MSE(t_s) = E(t_s - S_y^2)^2 = \left(S_y^2 e_0 - b_\varphi S_\varphi^2 e_1 - A_1 S_y^2 e_1\right)^2$$

$$= S_y^4 e_0^2 + \left[b_\varphi^2 S_\varphi^4 + A_1^2 S_y^4 + 2A_1 b_\varphi S_y^2 S_\varphi^2\right] e_1^2 - 2S_y^2 e_0 e_1 \left[b_\varphi S_\varphi^2 + A_1 S_y^2\right]$$

$$MSE(t_s) = \frac{1}{n}\left\{S_y^4(\lambda_{40} - 1) + (\lambda_{04} - 1)\left[b_\varphi^2 S_\varphi^4 + A_1^2 S_y^4 + 2A_1 b_\varphi S_y^2 S_\varphi^2\right] - 2S_y^2\left[b_\varphi S_\varphi^2 + A_1 S_y^2\right]\right\} \quad (3.7)$$

Table 3.1 presents some of the important estimators of the population variance, which can be obtained by suitable choice of constants $n_1, n_2$.

**Table 3.1: Members of $t_S$**

| Estimators | Values of | |
|---|---|---|
| | $n_1$ | $n_2$ |
| $t_{S1} = \dfrac{s_y^2 + b_\varphi(S_\varphi^2 - s_\varphi^2)}{s_\varphi^2} S_\varphi^2$ | 1 | 1 |
| $t_{S2} = \dfrac{s_y^2 + b_\varphi(S_\varphi^2 - s_\varphi^2)}{(s_\varphi^2 + \beta_{2\varphi})}(S_\varphi^2 + \beta_{2\varphi})$ | 1 | $\beta_{2\varphi}$ |
| $t_{S3} = \dfrac{s_y^2 + b_\varphi(S_\varphi^2 - s_\varphi^2)}{(s_\varphi^2 + C_p)}(S_\varphi^2 + C_p)$ | 1 | $C_p$ |
| $t_{S4} = \dfrac{s_y^2 + b_\varphi(S_\varphi^2 - s_\varphi^2)}{(s_\varphi^2 + \rho_{pb})}(S_\varphi^2 + \rho_{pb})$ | 1 | $\rho_{pb}$ |
| $t_{S5} = \dfrac{s_y^2 + b_\varphi(S_\varphi^2 - s_\varphi^2)}{(s_\varphi^2 \beta_{2\varphi} + C_p)}(S_\varphi^2 \beta_{2\varphi} + C_p)$ | $\beta_{2\varphi}$ | $C_p$ |
| $t_{S6} = \dfrac{s_y^2 + b_\varphi(S_\varphi^2 - s_\varphi^2)}{(s_\varphi^2 C_p + \beta_{2\varphi})}(S_\varphi^2 C_p + \beta_{2\varphi})$ | $C_p$ | $\beta_{2\varphi}$ |
| $t_{S7} = \dfrac{s_y^2 + b_\varphi(S_\varphi^2 - s_\varphi^2)}{(s_\varphi^2 C_p + \rho_{pb})}(S_\varphi^2 C_p + \rho_{pb})$ | $C_p$ | $\rho_{pb}$ |
| $t_{S8} = \dfrac{s_y^2 + b_\varphi(S_\varphi^2 - s_\varphi^2)}{(s_\varphi^2 \rho_{pb} + C_p)}(S_\varphi^2 \rho_{pb} + C_p)$ | $\rho_{pb}$ | $C_p$ |
| $t_{S9} = \dfrac{s_y^2 + b_\varphi(S_\varphi^2 - s_\varphi^2)}{(s_\varphi^2 \beta_{2\varphi} + \rho_{pb})}(S_\varphi^2 \beta_{2\varphi} + \rho_{pb})$ | $\beta_{2\varphi}$ | $\rho_{pb}$ |
| $t_{S10} = \dfrac{s_y^2 + b_\varphi(S_\varphi^2 - s_\varphi^2)}{(s_\varphi^2 \rho_{pb} + \beta_{2\varphi})}(S_\varphi^2 \rho_{pb} + \beta_{2\varphi})$ | $\rho_{pb}$ | $\beta_{2\varphi}$ |

Motivated by Singh et al. (2007), we propose another improved ratio- type estimator $t_{RS}$ for the population variance as

$$t_{RS} = s_y^2 \frac{(\eta S_\varphi^2 - v)}{[\alpha(\eta s_\varphi^2 - v) + (1-\alpha)(\eta S_\varphi^2 - v)]} \tag{3.8}$$

where $\eta, v$ are either real numbers or the functions of the known parameters of attributes such as $C_p, \beta_{2\varphi}, \rho_{pb}$ and $k_{pb}$.

Expressing equation (3.8) in terms of e's, we have

$$t_{RS} = S_y^2(1+e_0)\left[1 + A_2 e_1^2\right]^{-1}$$

where, $A_2 = \dfrac{\eta S_\varphi^2}{\eta S_\varphi^2 - v}$.

Up to first order of approximation, the MSE of $t_{RS}$ is given by

$$MSE(t_{RS}) = E(t_s - S_y^2)^2 = (S_y^2 e_0 - A_2 \alpha S_y^2 e_1)^2$$

$$MSE(t_{RS}) = \frac{S_y^4}{n}\left\{(\lambda_{40} - 1) + A_2^2 \alpha^2 (\lambda_{04} - 1) - 2A_2 \alpha(\lambda_{22} - 1)\right\} \tag{3.9}$$

Minimization of (3.9) with respect to $\alpha$ yields its optimum value as

$$\alpha_{opt} = \frac{(\lambda_{22} - 1)}{A_2(\lambda_{04} - 1)}.$$

Substituting optimum value of $\alpha$ in (3.9), we get the minimum variance of $t_{RS}$.

Table 3.2 presents some of the important estimators of the population variance, which can be obtained by suitable choice of constants $\eta, \nu$.

**Table 3.2: Members of $t_{RS}$**

| Estimator | Values of | | |
|---|---|---|---|
| | $\alpha$ | $\eta$ | $\nu$ |
| $t_{RS0} = s_y^2$ | 0 | 0 | 0 |
| $t_{RS1} = s_y^2 \dfrac{S_\varphi^2}{s_\varphi^2}$ | 1 | 1 | 0 |
| $t_{RS2} = s_y^2 \dfrac{(S_\varphi^2 - C_p)}{(s_\varphi^2 - C_p)}$ | 1 | 1 | $C_p$ |
| $t_{RS3} = s_y^2 \dfrac{(S_\varphi^2 - \beta_{2\varphi})}{(s_\varphi^2 - \beta_{2\varphi})}$ | 1 | 1 | $\beta_{2\varphi}$ |
| $t_{RS4} = s_y^2 \dfrac{(\beta_{2\varphi} S_\varphi^2 - C_p)}{(\beta_{2\varphi} s_\varphi^2 - C_p)}$ | 1 | $\beta_{2\varphi}$ | $C_p$ |
| $t_{RS5} = s_y^2 \dfrac{(C_p S_\varphi^2 - \beta_{2\varphi})}{(C_p s_\varphi^2 - \beta_{2\varphi})}$ | 1 | $C_p$ | $\beta_{2\varphi}$ |
| $t_{RS6} = s_y^2 \dfrac{(S_\varphi^2 + \beta_{2\varphi})}{(s_\varphi^2 + \beta_{2\varphi})}$ | 1 | 1 | $-\beta_{2\varphi}$ |

## 4. The suggested class of estimators

We suggest another improved class of estimators $t_M$ for population variance $s_y^2$ as

$$t_M = s_y^2 \left[ m_1 + m_2 \left( S_\varphi^2 - s_\varphi^2 \right) \right] \exp\left( \gamma \frac{[\delta S_\varphi^2 + \mu] - [\delta s_\varphi^2 + \mu]}{[\delta S_\varphi^2 + \mu] + [\delta s_\varphi^2 + \mu]} \right)$$

(4.1)

where $\delta$ and $\mu$ are either real numbers or function of known parameters of the auxiliary attribute $\varphi$ such as $C_p, \beta_{2\varphi}, \rho_{pb}$ and $k_{pb}$. The scalar $\gamma$ takes values -1 and +1 for ratio and product type estimators, respectively.

Expressing equation (4.1) in terms of e's and retaining terms up to second degree of e's, we have

$$t_M = S_y^2(1+e_0)[m_1 - m_2 s_\varphi^2 e_1]\left\{1 - \lambda\theta e_1 + \gamma\left(1+\frac{\gamma}{2}\right)\theta^2 e_1^2\right\} \tag{4.2}$$

Up to first order of approximation, the MSE of the estimator $t_M$ is

$$\text{MSE}(t_M) = S_y^4\left[1 + m_1^2 R_1 + m_2^2 R_2 + 2m_1 m_2 R_3 - 2m_1 R_4 - 2m_2 R_5\right] \tag{4.3}$$

where,

$$R_1 = 1 + \frac{1}{n}\left[(\lambda_{40}-1) + \gamma^2\theta^2(\lambda_{04}-1) + 2\gamma\left(1+\frac{\gamma}{2}\right)\theta^2(\lambda_{04}-1) - 4\gamma\theta(\lambda_{22}-1)\right],$$

$$R_2 = \frac{1}{n}S_\varphi^4(\lambda_{40}-1),$$

$$R_3 = \frac{1}{n}S_\varphi^2[2(\lambda_{22}-1) + 2\gamma\theta(\lambda_{04}-1)],$$

$$R_4 = 1 + \frac{1}{n}\left[\gamma\left(1+\frac{\gamma}{2}\right)\theta^2(\lambda_{04}-1) - \gamma\theta(\lambda_{22}-1)\right],$$

$$R_5 = \frac{1}{n}S_\varphi^2[\gamma\theta(\lambda_{04}-1) - (\lambda_{22}-1)].$$

On partially differentiating (4.3) with respect to $m_i (i=1,2)$, we get optimum values of $m_1$ and $m_2$, respectively as

$$m_{1(opt)} = \frac{(R_2 R_4 - R_3 R_5)}{R_1 R_2 - R_3^2}$$

and $\quad m_{2(opt)} = \dfrac{(R_1 R_5 - R_3 R_4)}{R_1 R_2 - R_3^2}.$

## 5. Efficiency Comparisons

First, we compare the efficiency of proposed estimator under optimum condition with usual estimator:

$$V(\hat{S}_y^2) - MSE(t_M) = \frac{S_y^4}{n}(\lambda_{40} - 1) - S_y^4\left[1 + m_1^2 R_1 + m_2^2 R_2 + 2m_1 m_2 R_3 - 2m_1 R_4 - 2m_2 R_5\right] \geq 0 \quad (5.1)$$

Next, we compare the efficiency of proposed estimator under optimum condition with the ratio estimator, exponential estimator, regression estimator and other estimators listed in the paper.

From (2.4) and (4.3), we have

$$MSE(t_1) - MSE(t_M) = \frac{S_y^4\left[(\lambda_{40} - 1) + (\lambda_{04} - 1) - 2(\lambda_{22} - 1)\right]}{n}$$

$$- S_y^4\left[1 + m_1^2 R_1 + m_2^2 R_2 + 2m_1 m_2 R_3 - 2m_1 R_4 - 2m_2 R_5\right] \geq 0 \quad (5.2)$$

From (2.5) and (4.3), we have

$$V(t_2) - MSE(t_M) = \frac{1}{n}\left[S_y^4(\lambda_{40} - 1) + b^2 S_\varphi^2(\lambda_{04} - 1) - 2b S_y^2 S_\varphi^2(\lambda_{22} - 1)\right]$$

$$- S_y^4\left[1 + m_1^2 R_1 + m_2^2 R_2 + 2m_1 m_2 R_3 - 2m_1 R_4 - 2m_2 R_5\right] \geq 0 \quad (5.3)$$

From (2.8) and (4.3), we have

$$MSE(t_3) - MSE(t_M) = \frac{S_y^4}{n}\left[(\lambda_{40} - 1) + \frac{(\lambda_{04} - 1)}{4} - (\lambda_{22} - 1)\right]$$

$$- S_y^4\left[1 + m_1^2 R_1 + m_2^2 R_2 + 2m_1 m_2 R_3 - 2m_1 R_4 - 2m_2 R_5\right] \geq 0 \quad (5.4)$$

From (3.5) and (4.3), we have

$$MSE(t_{KCi}) - MSE(t_M) = \frac{S_y^4}{n}\left[(\lambda_{40} - 1) + \omega_i^2(\lambda_{04} - 1) - 2\omega_i(\lambda_{22} - 1)\right]$$

$$- S_y^4\left[1 + m_1^2 R_1 + m_2^2 R_2 + 2m_1 m_2 R_3 - 2m_1 R_4 - 2m_2 R_5\right] \geq 0 \quad (5.5)$$

From (3.7) and (4.3), we have

$$MSE(t_s) - MSE(t_M) = \frac{1}{n}\left\{S_y^4(\lambda_{40} - 1) + (\lambda_{04} - 1)\left[b_\varphi^2 S_\varphi^4 + A_1^2 S_y^4 + 2A_1 b_\varphi S_y^2 S_\varphi^2\right] - 2S_y^2\left[b_\varphi S_\varphi^2 + A_1 S_y^2\right]\right\}$$

$$-S_y^4\left[1+m_1^2R_1+m_2^2R_2+2m_1m_2R_3-2m_1R_4-2m_2R_5\right]\geq 0 \quad (5.6)$$

$$\text{MSE}(t_{RS})-\text{MSE}(t_M)=\frac{S_y^4}{n}\left\{(\lambda_{40}-1)+A_2^2\alpha^2(\lambda_{04}-1)-2A_2\alpha(\lambda_{22}-1)\right\}$$

$$-S_y^4\left[1+m_1^2R_1+m_2^2R_2+2m_1m_2R_3-2m_1R_4-2m_2R_5\right]\geq 0 \quad (5.7)$$

Using (5.1)-(5.7), we conclude that the proposed estimator $t_M$ under optimum condition performs better than the other estimators discussed in this paper.

## 6. Empirical study

For empirical study, we use the data given in Sukhatme and Sukhatme (1970), p. 256.

The variables of the interest are:

Y is number of villages in the circle, and

$\phi$ represents a circle consisting more than five villages.

The values of required parameters are:

N = 89, n=23, $S_y^2$ = 4.074, $S_\phi^2$ =0.11, $C_y$= 0.601, $C_p$=2.678, $\rho_{pb}$ = 0.766,

$\beta_{2\phi}$ = 6.162, $\lambda_{22}$ = 3.996, $\lambda_{40}$ = 3.811, $\lambda_{04}$ = 6.162, $k_{pb}$ = 0.475931

**Table 6.1**: **PRE of various estimators**

| Estimator | PRE | Estimator | PRE |
|---|---|---|---|
| $S_y^2$ | 100 | $t_{S4}$ | 243.68 |
| $t_1$ | 141.89 | $t_{S5}$ | 218.49 |
| $t_2$ | 262.18 | $t_{S6}$ | 259.64 |
| $t_3$ | 254.27 | $t_{S7}$ | 191.17 |
| $t_{KC1}$ | 100.96 | $t_{S8}$ | 261.01 |
| $t_{KC2}$ | 99.58 | $t_{S9}$ | 126.42 |
| $t_{KC3}$ | 106.04 | $t_{S10}$ | 261.96 |
| $t_{KC4}$ | 101.10 | $t_{RS1}$ | 141.89 |
| $t_{KC2}$ | 99.58 | $t_{RS2}$ | 91.35 |
| $t_{KC3}$ | 106.04 | $t_{RS3}$ | 96.25 |
| $t_{KC4}$ | 101.10 | $t_{RS4}$ | 51.26 |
| $t_{S1}$ | 250.34 | $t_{RS5}$ | 90.06 |
| $t_{S2}$ | 261.80 | $t_{RS6}$ | 103.78 |
| $t_{S3}$ | 260.23 | $t_{RS}$ | **262.18** |

In Table 6.1, the percent relative efficiencies of the proposed estimators $t_{KCi}$ (i = 1,2,3,4), $t_{Si}$ (i = 1,2,...,10) and $t_{RSi}$ (i = 1,2,...,6) are listed when we choose different values of $n_1$ and $n_2$ in case of the estimator $t_{Si}$ and α, η and ν in case of the estimator $t_{RSi}$ respectively. Also, the PRE of regression type estimator $t_{RS}$ is 262.18.

**Table 6.2**: **PRE of suggested estimators $t_M$ with different values of constants**

| Values of δ and μ when (γ =1) | | |
|---|---|---|
| δ | μ | PRE |
| 1 | 1 | 264.54 |
| N | 1 | 274.03 |
| N | f | 264.57 |
| N | g=(1-f) | 274.17 |
| $β_{2φ}$ | $k_{pb}$ | 263.72 |
| 1 | 0 | **284.57** |
| N | $\rho_{pb}$ | 275.10 |
| N | $C_p$ | 278.75 |
| $\beta_{2\varphi}$ | $C_p$ | 263.18 |
| $C_p$ | $\beta_{2\varphi}$ | 265.43 |
| N | $k_{pb}$ | 270.93 |
| n | f | 278.65 |

In the Table 6.2, PRE of the proposed estimator $t_M$ with respect to $S_y^2$ is calculated for different values of parameters. It is observed the highest PRE 284.57 is obtained for $\delta = 1$ and $\mu = 0$. It has been also observed that the suggested class of estimators $t_M$ under optimum condition is more efficient than usual unbiased estimator, usual regression estimator, Singh and Kumar (2011) estimator and other estimators discussed in this paper. Hence, for proposed choice of parameters the proposed estimator $t_M$ is best among all the estimators considered in this paper.

**Acknowledgement**

The authors are very indebted to the editor in chief Prof. Melvin Scott and anonymous referee for their valuable suggestions leading to improvement of the quality of contents and presentation of the original manuscript.